\documentclass[12pt]{iopart}
\usepackage{iopams,amssymb,amsthm,latexsym,graphicx}
\usepackage{pictex}

\newtheorem{theorem}{Theorem}
\newtheorem{remark}[theorem]{Remark}
\newtheorem{definition}[theorem]{Definition}
\newtheorem{proposition}[theorem]{Proposition}
\newtheorem{corollary}[theorem]{Corollary}

\newcommand{\dom}{\mathop{\mathrm{dom}}\nolimits}
\newcommand{\ran}{\mathop{\mathrm{ran}}\nolimits}

\newcommand{\codim}{\mathop{\mathrm{codim}}\nolimits}
\newcommand{\ind}{\mathop{\mathrm{index}}\nolimits}
\newcommand{\G}{\Gamma}
\newcommand{\C}{\mathbf{C}}

\begin{document}

\title{Index theorems for quantum graphs}

\author[S Fulling, P Kuchment and J Wilson]
{S A Fulling$^{1,2,3}$, P Kuchment$^{1,3}$ and J H
Wilson$^{1,2,4}$}

\address{$^{1}$Department of Mathematics,
 Texas A\&M University, College Station, TX, 77843-3368, USA}
\address{$^{2}$Department of Physics, Texas A\&M University,
College Station, TX, 77843-4242, USA}
\address{$^{3}$Isaac Newton Institute for Mathematical Sciences,
University of Cambridge, 20~Clarkson Rd., Cambridge, CB3~0EH, UK
(during Spring, 2007)}
\address{$^{4}$Department of Physics, University of Maryland,
College Park, MD, 20742-4111, USA (from Fall, 2007)}

\eads{\mailto{fulling@math.tamu.edu}, \mailto{kuchment@math.tamu.edu},
\mailto{jwilson.thequark@gmail.com}}

\bigskip\centerline{24 August 2007}

\begin{abstract}
In geometric analysis, an index theorem  relates the difference
of the
numbers of solutions of two differential equations to the
topological structure of the manifold or bundle concerned, sometimes
using the heat kernels of two higher-order differential operators as
an intermediary. In this paper, the case of quantum graphs is
addressed. A quantum graph is a graph considered as a (singular)
one-dimensional variety and equipped with a second-order
differential Hamiltonian $H$ (a ``Laplacian'') with suitable
conditions at vertices. For the case of scale-invariant vertex
conditions (i.e., conditions that do not mix the values of functions
and of their derivatives), the constant term of the heat-kernel
expansion is shown to be proportional to the trace of the internal
scattering matrix of the graph. This observation is placed into the
index-theory context by factoring the Laplacian into two first-order
operators, $H =A^*A$, and relating the constant term to the index
of~$A$. An independent consideration provides an index formula for
any differential operator on a finite quantum graph in terms of the
vertex conditions. It is found also that the algebraic multiplicity
of $0$ as a root of the secular determinant of $H$ is the sum of the
nullities of
$A$ and $A^*$.
 \end{abstract}

\pacs{02.60.Lj, 02.40.-k, 02.30.Tb}
\ams{34B45, 47A53}
\submitto{\JPA}

\maketitle

\section{Introduction} \label{sec:intro}

A \emph{quantum graph} $\Gamma$ is a graph (with multiple edges
and loops allowed) in which each edge $e$ is assigned a coordinate
$x$ (and hence a length)
and the whole graph is equipped with a
self-adjoint differential operator as Hamiltonian. For instance, if
$\Gamma$ is embedded nicely into a Riemannian manifold, one can use
the arc length as a coordinate along an edge. In geometric language,
$\Gamma$ is a one-dimensional CW-complex with a Riemannian metric on
each 1-cell and appropriate boundary conditions at each 0-cell to
define a self-adjoint ``Laplacian''. We refer to \cite{GS, KSc, KSm,
Kuchment} for further background on quantum graphs.

Here we consider only graphs that are \emph{compact}: the number
of vertices $V$, the number of edges $E$, and the lengths of  all
edges are assumed to be finite. The number of edges attached to
vertex $v$ is called the \emph{degree} of~$v$ and denoted~$d_v\,$.
We also assume that every vertex has at least one edge attached,
since isolated vertices are negligible in the quantum graph
context.

The simplest Hamiltonian for a quantum graph is the
\emph{Laplacian} with \emph{Kirchhoff boundary conditions}, which
acts as the negative second derivative along each edge,
\begin{eqnarray*}
H=-\,\frac{\rmd^2}{\rmd x_e^2}\,,
\end{eqnarray*}
with the functions in its domain required to be continuous at
the vertices and to satisfy the Kirchhoff condition of no net flux
at each vertex:
\begin{eqnarray*}
\sum_{e\in E_v} {\rmd f\over \rmd x_e}(v) =0.
\end{eqnarray*}
Here $E_v$ is the set of edges incident on vertex $v$, and $x_e$
is the arc-length coordinate on $e$  outgoing from $v$
(in other words, the distance from $v$ of the variable point
on~$e$).
For more precise definitions see \cite{KSc, Kuchment}
and \sref{sec:graphs} below.

In one of the earliest papers on quantum graphs \cite{Roth},
J.-P.~Roth calculated the trace of the heat kernel $K$ for the
operator just described. He found an exact formula
\begin{equation}
 \sum_{n=0}^\infty \rme^{-\lambda_n t}
= \Tr K
= \int_\Gamma K(t,x,x)\, \rmd x
 = K_1 + K_2 + K_3\,,
\label{roth}\end{equation} where $\lambda_n$ are the eigenvalues of
$H$. Written in detail, it contains (before the integration) one
term for every path in the graph $\Gamma$ leading from the point $x$
(not a vertex) to itself. These closed paths fall into three
classes: The path of zero length yields the term $K_1 = 
L/\sqrt{4\pi
t}$, the anticipated leading term in the Weyl series, where $L$ is
the total length of all edges of the graph. $K_2$ is the sum of the
contributions of the periodic paths (where the initial and final
direction of the path are the same), which are proportional to
$\rme^{-L(C)^2/4t}$, where $L(C)$ is the length of the path $C$.
(Such terms do not contribute to the asymptotic expansion for $t\to
0$ of $\,\Tr K$ in powers of~$t$, but they determine oscillations in
the distribution of the eigenvalues~$\{\lambda_n\}$.) Finally, the
contributions of paths that are \emph{closed but not periodic}
(i.e., the initial and final directions are opposite) sum to the
simple form
\begin{equation}
K_3 = \case12 (V-E),
\label{Ktrace} \end{equation}
which constitutes the entire remainder of the Weyl series for the heat
kernel trace.
Note that $K_3$ is independent of $t$ and moreover is the \emph{only}
constant term in the formula~\eref{roth}.

The expression in \Eref{Ktrace} is interesting because it is a
half-integer and depends on the topology of $\Gamma$ only
(e.g., it is independent of the edge lengths).
 Indeed, it is the Euler
characteristic of $\Gamma$ regarded as a 1-complex. These features
are reminiscent of \emph{index theorems} in geometric analysis and
the calculation of indexes from the constant terms in heat-kernel
expansions \cite{Gilkeybook,Gilkey}. The original goal of this
paper was to give an index interpretation of~\eref{Ktrace}; in
fact, we also generalize it to graphs with other boundary
conditions and  compute indexes of quite general quantum graph
operators by another, very simple method.

The contents of the paper are as follows:
\Sref{sec:interval} reviews the
appearance of indexes in the heat kernel asymptotics for the case
of an interval. In  \sref{sec:graphs} we introduce necessary
notions and auxiliary results concerning quantum graphs. The
section also contains a general formula for indexes of
differential operators on quantum graphs. \Sref{sec:index}
contains the main results concerning the relations between the
constant terms in the asymptotic expansion of the heat kernel and
indexes of suitable operators on the graph.
 Relations to the secular determinant are
discussed in  \sref{sec:secular}. The final
\sref{sec:concl} contains some remarks and conclusions.

\section{The interval} \label{sec:interval}

In \cite[Section 1.5]{Gilkey} P.~Gilkey treats the Laplacian on an
interval with Dirichlet and with Neumann boundary conditions as
the prototype of the index theorem for the de~Rham complex on a
manifold with  boundary. The index theorem for a quantum graph
with Kirchhoff boundary conditions is a different generalization
of this elementary example, so we shall review the latter.

Let $H_D$ and $H_N$ be the operator $- \rmd^2/\rmd x^2$ on the
interval $(0,L)$ with Dirichlet and Neumann boundary conditions,
respectively. The eigenfunctions of $H_N$ with eigenvalue $0$ are
the constant functions, so the kernel of $H_N$ has dimension~$1$. In
contrast, $H_D$ has trivial kernel. On the other hand
\cite[Subsection 3.1.3]{Gilkey}, the heat traces of these operators
are
\begin{equation}
\Tr\,\rme^{-tH_{D,N}}  = \int_0^L K_{D,N}(t,x,x)\,\rmd x
\mathop{\sim}\limits_{t\to 0} {L\over\sqrt{4\pi t}} \mp \frac12\,,
\label{intvaltrace}\end{equation} where the negative sign applies to
the Dirichlet case and the exponentially small terms analogous to
$K_2$ in \eref{roth} have been omitted. Therefore,
\begin{equation}
\dim\ker H_N - \dim\ker H_D = 1 = \Tr K_N - \Tr K_D\,.
\label{NDindex}\end{equation}

To identify \eref{NDindex} as an index theorem we must factor
$H_{D,N}$ into first-order operators. Let $A$ be the operator
$\rmd/\rmd x$ acting on the domain $H^1(0,L)$, which is the Sobolev
space containing functions on $(0,L)$ that, together with their
first distributional derivatives, are square-integrable.
Standard
integrations by parts show that the adjoint operator $A^*$ is
$-\rmd/\rmd x$ with domain $H^1_0(0,L)$ (containing functions from
$H^1(0,L)$ that satisfy the Dirichlet conditions $f(0)=f(L)=0$) and
that $A^{**}=A$. One now forms two second order operators
\begin{equation}
H_N=A^*A, \qquad H_D = AA^*, \label{intval2ndorder}\end{equation}
where in the first case the domain consists of twice differentiable
($H^2(0,L)$) functions satisfying the Neumann conditions,
$f'(0)=f'(L)=0$, so that $Af\in \dom A^*$ and the composition is
defined; in the second case, similarly, the domain consists of
functions from
$H^2(0,L)$ satisfying the Dirichlet conditions.

Because $\langle f,A^*Af\rangle = \langle Af,Af\rangle = \|Af\|^2$,
the kernel of $A$ is the same as that of $H_N\,$.
Similarly, $\,\ker A^* =\ker H_D\,$.
Therefore, \eref{NDindex} can be restated
as the index formula
\begin{equation}
\ind A = \Tr K_N - \Tr K_D = 1.
\label{intvalindex}\end{equation}

\begin{remark} \label{oneforms}\rm
By regarding the elements of $\,\dom A^*$ as 1-forms rather than
scalar functions, one identifies $A$ and $A^*$ with the exterior
derivative operator $\,\rmd\colon \Lambda^0(0,L)\to\Lambda^1(0,L)$
and its adjoint $\delta\colon \Lambda^1(0,L)\to\Lambda^0(0,L)$,
where $\Lambda^0(0,L)$ and $\Lambda^1(0,L)$ are the spaces of
$L^2$ functions $f(x)$ and $1$-forms $g(x)\,dx$. One can combine
these
operators into a single operator from $\Lambda^0(0,L) \oplus
\Lambda^1(0,L)$ into itself,
\begin{eqnarray*}
\rmd+\delta =\left(\begin{array}{cc}
 0 & \delta \\ \rmd & 0
\end{array}\right),
\end{eqnarray*}
where, in the version on the left, $\,\rmd\,$ annihilates the 1-forms and
$\delta$ annihilates functions.
 Then the Hodge Laplacian is
\begin{eqnarray*}
(\rmd+\delta)^2 = \delta\rmd+\rmd\delta
 = \left(\begin{array}{cc}
 H_N &0 \\ 0& H_D
\end{array}\right).
\end{eqnarray*}
It is this formulation that generalizes to higher-dimensional manifolds,
with $H_N$ acting on forms of even degree and $H_D$ on forms of odd
degree (or vice versa)
\cite{Gilkeybook,Gilkey}.
\end{remark}

In the following sections
we will  extend this analysis to more general quantum graphs.
In particular, a central task is to identify the analogues of the
operators $A$ and~$A^*$.

\section{Quantum graphs}\label{sec:graphs}

\subsection{Vertex conditions} \label{ssec:vertex}
As we have mentioned in \sref{sec:intro}, appropriate vertex
conditions are needed in order to turn the (negative) second
derivative along the edges of a quantum graph into a self-adjoint
operator in $L^2(\G)$. All such choices of boundary conditions at
vertices were catalogued in
\cite{KSc} (after prior discussion in \cite{Exner}) and
reformulated in \cite{Har, Kuchment}.
It will be convenient for us to follow the
formulation from \cite{Kuchment}.

Let $v$ be a vertex and $f(x)$ a
function on $\G$. We denote by
\begin{equation*}
F(v) = \left(\begin{array}{c} f_1 (v)\\ \vdots \\
f_{d_v}(v)
\end{array}\right)
\end{equation*}
the vector of values of the function $f$ at the vertex $v$,
attained along $d_v$ edges incident to $v$. In particular, if $f$
were continuous, all these values would be equal.
Analogously,
\begin{equation*}
F'(v) = \left(\begin{array}{c} f'_1(v) \\ \vdots \\
f'_{d_v}(v)
\end{array}\right)
\end{equation*}
is the vector of derivatives at $v$ of $f$ along these edges,
where the derivatives are taken in the directions outgoing from
the vertex $v$.

It is clear that vertex conditions for the second-derivative
operator can involve only the values of the function and of its
derivatives along edges. If these conditions do not mix the values
attained at different vertices, they are called \emph{local}.
(On an interval, for instance, Dirichlet,
Neumann, and Robin conditions are local, but the periodicity
condition is nonlocal.)
 As we will see later, there is actually not much
difference between local and nonlocal vertex conditions on a
quantum graph. (For instance, the periodicity condition becomes
local if the interval is replaced by a loop attached to a single
vertex.)

\begin{theorem} {\rm\cite{Kuchment}} \label{T:Sa}
All self-adjoint realizations $H$ of the negative second
derivative on $\Gamma$ with local vertex boundary conditions can
be described as follows: For every vertex $v$, of degree~$d_v\,$,
there are two orthogonal (and mutually orthogonal) projectors
$P_v$, $Q_v$ operating in $\mathbb{C}^{d_v}$ and an invertible
self-adjoint operator $\Lambda_v$ operating in the subspace
$(1-P_v-Q_v)\mathbb{C}^{d_v}$. (Either $P_v$, $Q_v$, or
$C_v \equiv 1-P_v-Q_v$ might be zero.) The functions $f$ in the
operator
domain are those members of the Sobolev space $\bigoplus_e H^2(e)$
that satisfy at each vertex $v$ boundary conditions consisting of
the ``Dirichlet part''
\begin{equation} P_vF(v)=0, \label{dirpart}
\end{equation}
the ``Neumann part''
\begin{equation} Q_vF'(v)=0,
\label{neupart}
\end{equation}
and the ``Robin part''
\begin{equation}
C_vF'(v) =\Lambda_v C_vF(v).\label{robpart}
\end{equation}

The quadratic form of $H$ is
\begin{equation}\label{qform}
 h[f,f]=\sum\limits_{e\in E}\int\limits_e
\left|\frac{df}{dx}\right|^2\,dx
+\sum\limits_{v\in V}(\Lambda_v C_v F,C_v F) \end{equation}
with the domain that consists of the functions $f(x)$ that belong
to the Sobolev space $H^1(e)$ on each edge $e$ and satisfy
\eref{dirpart} at each vertex.
\end{theorem}
\begin{remark}{\em
This theorem was formulated a little bit
differently  in \cite{Kuchment}.
For one thing, $\Lambda_v$ was called $-L_v$ there.
More importantly, there the two
projectors $Q_v$ and $C_v$ where lumped into a single one, and
thus the condition of invertibility of the operator $\Lambda_v$
disappeared.
 The equivalent formulation provided here
 distinguishes between Robin and pure Neumann conditions,
as is often convenient.}
\end{remark}
\begin{remark}
{\em The three parts (\ref{dirpart})--(\ref{robpart}) of the
vertex conditions can be combined into a single condition
\begin{equation}\label{E:condit}
    A_vF(v)+B_vF^\prime (v)=0,
\end{equation}
where the $d_v\times d_v$ matrices $A_v$ and $B_v$ are
\begin{equation}\label{E:matrix}
A_v=P_v-\Lambda_vC_v, \quad  B_v= Q_v+C_v\,.
\end{equation}
The conditions were introduced in \cite{KSc}
in the form (\ref{E:condit})
(which by itself does not suffice to define the
matrices $A_v$ and $B_v$  uniquely, however).}
\end{remark}

The most popular  vertex conditions are the \emph{Kirchhoff}
ones (also called \emph{Neumann} or \emph{natural}),
which reduce at vertices of degree $1$  to
Neumann conditions:
\begin{definition} \label{def:kirch}
The \emph{Kirchhoff} boundary conditions are defined by \numparts
the continuity condition
\begin{equation} f_1(v) = f_2(v) = \cdots
= f_{d_v}(v) \equiv f(v) \label{kirchdir}\end{equation} as
Dirichlet part and the zero flux condition
\begin{equation} \sum_{e=1}^{d_v} f'_e(v) = 0
\label{kirchneu}\end{equation} as Neumann part, with no Robin
part.
\endnumparts
\end{definition}

In other words, $(1-P_v)\mathbb{C}^{d_v}$ is in this case
one-dimensional and consists of  the vectors with equal
coordinates.

Another type of conditions that arises in our work
is dual to the Kirchhoff type,
in the sense that  the roles
of the values and derivatives of the function $f$ at each vertex
are switched.
(At vertices of degree $1$ these ``anti-Kirchhoff'' conditions
reduce to the Dirichlet ones.)

\begin{definition} \label{def:anti}
The \emph{anti-Kirchhoff} boundary conditions are \numparts
\begin{equation} \sum_{e=1}^{d_v} f_e(v) = 0
\label{antidir}\end{equation} as Dirichlet part and
\begin{equation}
f'_1(v) = f'_2(v) = \cdots = f'_{d_v}(v)\equiv f'(v)
\label{antineu}\end{equation} as Neumann part, with no Robin part.
\endnumparts
\end{definition}

\subsection{Bonds vs edges} \label{ssec:bonds}
In what follows, we will need to use directed edges (which we will
call {\em bonds}) rather than the undirected ones as before. Thus,
each edge results in {\em two} directed bonds (with opposite
directions), which are denoted by Greek letters. We  denote
by $\overline{\alpha}$ the bond $\alpha$ with its direction
reversed.

 Recall that loops (tadpoles) can always be removed from a quantum
graph by inserting extra Kirchhoff vertices of degree~$2$.
 Adding such a vertex does not change the heat trace or the Euler
characteristic, nor either side of any of the index formulas in this
paper.
 Therefore, one may assume that the two ends of a bond are distinct
vertices.

It is not necessary to pick either of the two directions of an edge
as the canonical one.  The language of differential forms makes it
possible to give global meaning to the differential of a function on
$\Gamma$ without committing to any particular coordinate, $x_e\,$,
on each edge.  In discussing the behavior of functions (and their
derivatives) in the neighborhood of any one vertex, therefore, we
remain free to use the most convenient coordinate on each edge,
namely, the outgoing arc length parameter.

\subsection{Scattering matrices}\label{ssec:scattering}
In this subsection we  introduce, following \cite{KSc,KSm},
the
scattering matrices and some of their properties that we will need
in the rest of the text.

Let $H$ be a self-adjoint realization of the negative second
derivative $-\rmd^2/\rmd x^2$ on a finite quantum graph~$\G$
(i.e., one of
the self-adjoint vertex conditions described in Theorem
\ref{T:Sa} is imposed).

Let us consider a vertex $v$ and the set $E_v$ of all edges $e$
incident to it.
 (Such a configuration is called a {\em star}; see
Fig.~\ref{F:star}.)
\begin{figure}[ht]
\[ \beginpicture \setcoordinatesystem units <.2in,.2in>
 \setplotarea x from -1.5 to 2, y from -1.2 to 1.5
 \put{$\bullet$} at 0 0
 \plot 0 0   2 0 /
 \plot 0 0   2 1.5 /
\plot 0 0   -1 1.5  /
\plot 0 0   -1.5 -1.2 /
\plot 0 0   1.3 -1 /
\endpicture \]
\caption{A star.}\label{F:star}
\end{figure}
For any edge $e_0\in E_v$ and any real $k$, we choose as in
\sref{sec:intro} the
coordinate $x$ increasing away from the vertex and consider the
unique solution $f(x)$ on the star $E_v$ of the following
scattering problem at $v$:
\begin{equation}\label{E:scatter}
    \cases{-\, \frac{\rmd^2f}{\rmd x^2}=k^2f(x) \mbox{ on each
edge }e\in
    E_v\,,
    \cr
f(x)=e^{-ikx}+\sigma^{(v)}_{e_0e_0}e^{ikx} \mbox{ on } e_0\,, \cr
f(x)=\sigma^{(v)}_{e_0e}e^{ikx} \mbox{ on } e\neq e_0\,, \cr
\mbox{vertex conditions are satisfied at }v.\cr}
\end{equation}
In other words, $\sigma^{(v)}_{e_0e_0}$ is the reflection
coefficient along the bond $e_0\,$, and $\sigma^{(v)}_{e_0e}$ is
the
transmission coefficient from the edge $e_0$ to $e$. Notice that
the coefficients $\sigma^{(v)}_{e_0e}$ in general depend on $k$.
\begin{definition}
The unitary $d_v\times d_v$ matrix $\sigma^{(v)}(k)$ with
the entries $\sigma^{(v)}_{e_1e_2}(k)$ for $e_j\in E_v$ is the
{\em edge scattering matrix at the vertex $v$}.
\end{definition}
Notice that in defining $\sigma^{(v)}$
the direction chosen along each edge
depends on the vertex considered.
That is why it becomes necessary to deal with directed bonds when
a scattering matrix $S$ for the whole graph is defined.
However, as explained in \sref{ssec:bonds}, the ambiguity in
$x_e$ does not create any inconsistency in the notation.
 Another remark is that this
matrix clearly depends upon what type of vertex conditions are
imposed,
and not every matrix function $\sigma^{(v)}(k)$ can necessarily
be realized by one of the second-order differential Hamiltonians
studied here.

It is straightforward to derive the  formula \cite{KSc}
\begin{equation}\label{E:scat}
    \sigma^{(v)}(k)=-(A_v+ikB_v)^{-1}(A_v-ikB_v),
\end{equation}
which in particular confirms that the matrix is unitary
and shows that its $k$-dependence is tightly constrained.
From (\ref{E:matrix}) we get an alternative representation of
$\sigma$:
\begin{equation}\label{E:scat_projectors}
\sigma^{(v)}(k)=P_v-Q_v+(\Lambda_v-ik)^{-1}(\Lambda_v+ik)C_v\,.
\end{equation}

The following result will be important for what follows. A part of
it was proved by Kostrykin and Schrader (\cite[Proposition
2.4]{KPS}), \cite[Corollary 2.3]{KSc}, \cite[Theorem 1]{KSc2}).

\begin{theorem}\label{T:ssquared}
The following conditions are equivalent:
\begin{enumerate}
\item For each vertex $v$, $\sigma^{(v)}$ is independent of $k$.

\item For each vertex $v$, there is a value $k\neq 0$ such that
$(\sigma^{(v)})^2 = 1$.

\item For each vertex $v$,   $(\sigma^{(v)})^2  = 1$ for all
$k$.

\item For each vertex $v$,  $\sigma^{(v)}$ has the form $1-2Q_v$
for some orthogonal projection $Q_v\,$.

\item There is no Robin part in the vertex conditions:
$C_v=0$ for each vertex $v$.

\item The vertex conditions are \emph{scale-invariant}
(i.e., if a
function $f(x)$ on
 neighborhood of $v$ in the star $E_v$
 satisfies the vertex conditions at~$v$,
 then after rescaling to $f(rx)$, it still satisfies the
conditions).

\item The Hamiltonian $H=-\frac{d^2}{dx^2}$ with the given vertex
conditions can be factored as $H=A^*A$, where $A=\frac{d}{dx}$
with appropriate vertex conditions, and $A^*$ is its adjoint
operator.
\end{enumerate}
\end{theorem}
\begin{proof} Equivalence of statements (i) through (iv) was
proved by Kostrykin and Schrader; it also follows rather easily
from \eref{E:scat_projectors}.

Computing $(\sigma^{(v)})^2$ using (\ref{E:scat_projectors}), we
get
\begin{equation}\label{E:ssquared}
(\sigma^{(v)})^2=P_v+Q_v+\left(\frac{\Lambda_v+ik}
{\Lambda_v-ik}\right)^2C_v\,.
\end{equation}
If $C_v\neq 0$, then in order to get $\sigma^2=1$, we need that
$\left(\frac{\Lambda_v+ik}{\Lambda_v-ik}\right)^2=1$. A
straightforward calculation shows that this is impossible for an
invertible operator $\Lambda$ and a non-zero $k$. This proves
equivalence of (iii) and (v).

Equivalence of (v) and (vi) is trivial.

The implication (vii) $\Rightarrow$ (v) can be established as
follows. The vertex conditions for $A$ can
involve only the values of the function, not its derivatives.
Thus, at any vertex $v$
they can be written as $P_v F(v)=0$ for some orthogonal projector
$P_v\,$. Then a simple and well known (e.g., \cite{KSc})
calculation, which boils down to an integration by parts, shows
that for the adjoint operator $-\frac{\rmd}{\rmd x}$\,, the
vertex
conditions are given by the complementary projector
$Q_v=1-P_v\,$.
The equality $H=A^*A$ now implies that the vertex conditions for
$H$ have $P_v$ as the Dirichlet and $Q_v$ as the Neumann part,
with no Robin part being present. This argument can easily be
reversed to show the converse implication, (v) $\Rightarrow$
(vii).
\end{proof}
\begin{corollary}
The scattering matrices $\sigma$ for the Hamiltonian
$-\frac{\rmd^2}{\rmd x^2}$ with Kirchhoff or anti-Kirchhoff
boundary
conditions satisfy the equivalent conditions of Theorem
\ref{T:ssquared}.
\end{corollary}
We now introduce the  \emph{global scattering matrix}
$S_{\alpha\beta}\,$, entries of which are indexed by the
directed bonds $\alpha$ and $\beta$.
\begin{definition}\label{D:scat}
  The $2E\times 2E$ \emph{(global) scattering matrix} $S$ is
defined as follows:
\begin{equation}\label{E:def_scat}
    S_{\beta\alpha}=\cases{ \sigma^{(v)}_{\beta\alpha} \mbox{ if $\alpha$
    terminates at $v$ and $\beta$ starts at
    $v$},
    \cr
    0 \mbox{ otherwise},
    }
\end{equation}
where $\alpha$ and $\beta$ are (directed) bonds in $\G$.
\end{definition}

\begin{proposition}\label{L:odd}
Let $H$ satisfy the equivalent conditions of Theorem
\ref{T:ssquared}, so that it factors 
as $H=A^*A$ as in (vii) of
the theorem. Let also $H^\prime\equiv AA^*$. We also denote by
$\sigma^\prime$ and $S^\prime$ the scattering matrices for
$H^\prime$. Then
\begin{equation}\label{E:odd}
\begin{array}{c}
   
\sigma^{(v)}_{\alpha\beta}=-(\sigma^\prime)^{(v)}_{\alpha\beta}\,,\\
S_{\alpha\beta}=-(S^\prime)_{\alpha\beta}
\end{array}
\end{equation}
for any vertex $v$ and any bonds $\alpha,\beta$.
\end{proposition}
\begin{proof}
Indeed, it is clear that for $H^\prime$ the projectors $P_v$ and
$Q_v$ exchange their places, while $C_v=0$. Then formulas
(\ref{E:scat_projectors}) and (\ref{E:def_scat}) prove the
statement.
\end{proof}
\subsection{Indexes of quantum graph operators}\label{sec:count}

As it happens, one can establish a simple formula for the index of
any (elliptic) differential operator on a compact quantum graph
$\Gamma$, which in particular implies the index formulas for the
exterior-derivative operators $A$ introduced previously.

First we need to review the basic notions
concerning the Fredholm property and the index (see, e.g.,
\cite{Kato}).
Recall that the codimension of a (closed) subspace
$E\subset H$ is defined as the dimension of the quotient space
$H/E$, or, equivalently (in a Hilbert space), the dimension of an
orthogonal complement of~$E$.

\begin{definition}
A bounded operator $T\colon H_1\to H_2$ between two Hilbert (or
Banach) spaces is said to be {\em Fredholm}, if it has a closed
range and the dimension of its kernel $\ker T$ and the codimension
of its range $\ran T$ are finite.
The \emph{index} of a Fredholm operator $T$ is defined as
\[
\ind T= \dim\ker T - \codim\ran T.
\]
\end{definition}

\begin{proposition}
If operator $T$ is Fredholm and operator $K$ is compact (in
particular, of  finite rank), then $T+K$ is also Fredholm and
$\ind (T+K)=\ind T$.
\end{proposition}

 To formulate the main theorem for a differential
operator of arbitrary order it is convenient to choose an
orientation for each edge, so that the arc length parameter $x_e$
is unambiguous.

\begin{theorem} \label{genorder} Consider the operator on $\Gamma$
defined by the differential expression of order~$m$
\begin{equation}
T = \sum_{j=0}^m c_j(x_e) \,{\rmd^{m-j}\over \rmd x_e^{m-j}}
\label{arborder}\end{equation} with $c_0(x)$ continuous on each
closed edge (but not necessarily on the whole graph) and never
equal to~$0$ and all other $c_j$ measurable and bounded. Let $T_1$
be the restriction of $T$ as an operator from $\oplus_e H^m(e)$
into $L^2(\Gamma$) to a subspace of codimension $p$ (e.g., by
imposing $p$ vertex conditions sustainable by $H^m$, i.e.,
involving derivatives up to the order $m-1$). Then
\begin{enumerate}
\item The operator so defined  is Fredholm.

\item
\begin{equation}
\ind T_1 = mE -p.
\end{equation}
\end{enumerate}
\end{theorem}
\begin{proof}
Consider $T$ as the naturally defined (and obviously bounded)
operator from $\oplus_e H^m(e)$ into $L^2(\Gamma)$.
 All terms
in $T$ that involve derivatives of order less than $m$ are compact
operators and thus do not influence the Fredholm property or the
index. Thus, we can assume that $T=c_0(x_e)\frac{\rmd^m}{\rmd x_e^m}$.
This is now the composition of $\frac{\rmd^m}{\rmd x_e^m}$ acting from
$\oplus_e H^m(e)$ to $L^2(\Gamma)$ with the invertible operator of
multiplication by $c_0(x)$ in $L^2(\Gamma)$. Thus, everything
reduces to the $m$th derivative alone. It is easy to show that it
is a surjective operator from $H^m(e)$ onto the whole $L^2(e)$
(and thus from  $\oplus_e H^m(e)$ onto $L^2(\Gamma)$). On each
edge, it has the $m$-dimensional kernel consisting of polynomials
of degree less than $m$. Thus, $T$ is Fredholm and
\[
\ind T=\sum_e (m - 0)=mE.
\]

Let us now notice that by definition, $T_1$ is the restriction of
$T$ onto a subspace $M$ of codimension $p$. Consider any
($p$-dimensional) complement $N$ to $M$ in $\oplus_e H^m(e)$ and
the extension $\tilde{T}$ of $T_1$ from $M$ to the whole $\oplus_e
H^m(e)$ that acts as the zero operator on $N$. Then the difference
$T-\tilde{T}$ vanishes on $M$ and thus is a finite-dimensional
operator. Hence, $\tilde{T}$ is Fredholm of the same index $mE$ as
$T$.

On the other hand, it is clear that the ranges of $\tilde{T}$ and
$T_1$ are the same and the kernel of $\tilde{T}$ is $p$ dimensions
larger than the kernel of $T$. Hence, $T$ is Fredholm and $\ind
T=\ind \tilde{T}-p=mE-p$.
\end{proof}
This implies in particular
\begin{corollary}\label{C:index1}
Let $A$ be the exterior derivative $\rmd$ acting as a bounded
operator from a subspace $M$ of codimension $p$ in $\bigoplus_e
H^1(e)$ into $L^2(\Gamma)$.  This operator is Fredholm of index
$E-p$.  Therefore,
\begin{equation}
 \ind A=E-p=E-\sum_v\dim P_v\,,
 \label{projindthm} \end{equation}
  where $P_v$ are the
orthogonal projectors describing the vertex conditions for $A$.
Thus, in particular,
 \begin{enumerate}

\item Without any vertex conditions, one has $\,\ind A =E$.

\item With continuity conditions \eref{kirchdir} at all vertices,
one has
 \[ \ind A =E-\sum_v(d_v-1)=E-(2E-V)=V-E. \]

\item With the condition \eref{antidir} that the sum of values at
each vertex is equal to zero, one has
 \[\ind A =E-\sum_v 1=E-V.\]
\end{enumerate}
\end{corollary}

 Notice that $\sum_v\dim P_v$ arising in this corollary is just the
number of vertex conditions defining $H$ that contain only the
values of the function and no derivatives (Dirichlet part of the
conditions).

\section{Heat kernel and index}\label{sec:index}
In this section we address the
relation between the heat trace asymptotics and the index on
quantum graphs. Most of the considerations are independent of
 Theorem~\ref{genorder}.

 We will assume from now on that $\G$ is an
arbitrary quantum graph and the Hamiltonian $H$ satisfies the
conditions of the Theorem \ref{T:ssquared}, so that it
factors as $H=A^*A$, where $A=\rmd/\rmd x$
  with some vertex conditions on
the values of functions, such conditions corresponding at any
vertex $v$ to an orthogonal projector $P_v$ in $\C^{d_v}$. Then,
as before, we denote by $H^\prime$ the operator $AA^*$ with the
vertex conditions given by the orthogonal projector 
$Q_v=1-P_v\,$.
Let also $K$ and $K^\prime$ be the corresponding heat kernels.
 The following proposition is standard.

\begin{proposition} \label{H_ind}
In the situation just described,
\begin{equation}
\ind A = \Tr K - \Tr K^\prime.
 \label{E:H_ind}
\end{equation}
\end{proposition}

\begin{proof}
The non-zero eigenvalues of $H=A^*A$ and $H^\prime=AA^*$  are the
same, including their multiplicity. The only exception is that the
dimensions of the eigenspaces for the eigenvalue~$0$ are
different. Thus, the difference of the heat kernel traces is
guaranteed to be independent of $t$. At large $t$ this difference
reduces to the difference of the nullities (i.e., dimensions of
the kernels), and at small $t$ it reduces to the difference of the
constant terms in the heat-kernel expansions. See, for instance,
\cite{Gilkeybook} for a more detailed exposition.
Now, since the nullity of $H$ is clearly equal to that of $A$, and
that of $H^\prime$ coincides with that of $A^*$, one concludes
that
$
\Tr K - \Tr K^\prime=\dim \ker H - \dim \ker H^\prime
=\dim \ker A - \dim \ker A^*=\ind A.
$
 \end{proof}

We now need to establish a formula for the constant term of a heat
trace:

\begin{theorem} {\rm\cite{Wilson,BHW}} \label{Ssum}
Let $\Gamma$ be a finite quantum graph with scale-invariant vertex
conditions defining the Laplacian.
  (Thus the bond-to-bond scattering matrix
$S_{\alpha\beta}$ is independent of the frequency $k$).
Let $K(t,x,y)$ be the corresponding heat kernel on~$\G$.
 Then the constant term
in the asymptotic expansion at $t\to 0$ of the heat trace
\begin{eqnarray}
 \sum_{n=0}^\infty \rme^{-\lambda_n t}
= \int_\Gamma K(t,x,x)\, \rmd x \label{trace}\end{eqnarray} is
\begin{equation}
\frac14 \sum_{\alpha}S_{\alpha\overline{\alpha}}\,.
\label{strace}
\end{equation}
(See Definition \ref{D:scat} for $S_{\alpha\beta}$.)
\end{theorem}

\begin{proof}[Sketch of proof]
 The theorem is proved in
\cite{BHW} for a different kernel, but as stressed in
\cite{Wilson} the same argument
applies to a whole class of kernels, including the heat kernel.
(See also \cite{Roth,KPS}.)
Starting from  the standard one-dimensional
heat kernel
 on the real line,
\begin{equation}\label{E:heat}
K_0(t,x,0) \equiv (4\pi t)^{-1/2} e^{-x^2/4t},
\end{equation}
by an extension of the method of images one constructs the heat
kernel on the graph as a sum over all paths from $y$ to~$x$, which
then needs to be restricted to the diagonal $y=x$. The heat trace
is formed then by integrating over $x$. As in \cite{Roth}, the
contributions of the periodic paths, i.e., the ones that return
to
the point $x$ with the same direction as at the start, are
proportional to Gaussian terms $e^{-L_\mathbf{q}^2/4t}$ and cannot
contribute to the $t$-independent term of the heat-kernel
expansion. The path of zero length gives the leading Weyl term,
proportional to $t^{-1/2}$.
\begin{figure}
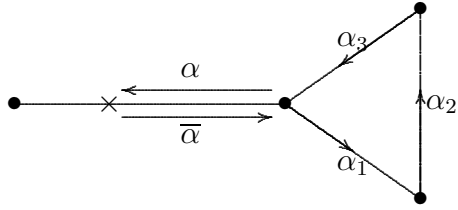

\[
\beginpicture
\setcoordinatesystem units <1.8truecm,1.8truecm> point at 0 0
\put{$\bullet$} at -2 0 \put{$\bullet$} at 0 0 \put{$\bullet$} at
1 0.7 \put{$\bullet$} at 1 -0.7 \put{$\times$} at -1.3 0 \plot -2
0
       0 0 /
\plot  1 0.7
       0 0 /
\plot  1 -0.7
       0 0 /
\plot  1 0.7
       1 -0.7 /
\arrow <5pt> [.2,.67] from 0 0 to 0.5 -0.35 \put{$\alpha_1$} [t]
<0pt,-2pt> at 0.5 -0.35 \arrow <5pt> [.2,.67] from 1 0.7 to 0.4
0.28 \put{$\alpha_3$} [t] <0pt,9pt> at 0.5 0.35 \arrow <5pt>
[.2,.67] from 1 -0.7 to 1 0.1 \put{$\alpha_2$} [l] <2pt, 0pt> at 1
0 \arrow <5pt> [.2,.67] from -0.1 0.1 to -1.2 0.1 \put{$\alpha$}
[b] <0pt,5pt> at -0.7 0.1 \arrow <5pt> [.2,.67] from -1.2 -0.1 to
-0.1 -0.1 \put{$\overline\alpha$} [t] <0pt,-2pt> at -0.7 -0.1
\endpicture \]
\caption{A bounce path
$\alpha\alpha_3\alpha_2\alpha_1\overline\alpha$ with $n=3$.}
\label{fig:stirrups}
\end{figure}
The contribution of the other class of paths, which are closed but
not periodic (``bounce'' paths),
$\alpha\mathbf{p}\overline\alpha$, where $\mathbf{p}$ is a cycle
in~$\G$ (see Fig.~\ref{fig:stirrups}, where the triangle represents
the cycle $\mathbf{p}$ and the point $x$ is located on the
bond~$\alpha$),
can be reduced to the following sum:
\begin{equation}\label{E:sum}
    \frac{1}{2}\sum_n\sum_{\mathbf{p}\in P_n}\sum_\alpha
A_{\alpha\mathbf{p}\overline\alpha}
\int^{l_\mathbf{p}+2L_\alpha}_{l_\mathbf{p}}K_0(t,x)\,dx.
\end{equation}
(The condition of $k$-independence of the scattering matrix is
used
here.) We denote here by $P_n$ the set of cycles of period
(number
of edges traversed, including multiplicity)~$n$.
E.g., in Fig.~\ref{fig:stirrups},
 the cycle $\mathbf{p}=\alpha_1\alpha_2\alpha_3$ has
period $3$. We also use here notations $L_\alpha$ for the length
of the bond $\alpha$ and $l_\mathbf{p}$ for the metric length of
the cycle
$\mathbf{p}$ (e.g.,
$l_\mathbf{p}=L_{\alpha_1}+L_{\alpha_2}+L_{\alpha_3}$
in
 Fig.~\ref{fig:stirrups}). The shorthand notation $A_{\alpha
\mathbf{p}\overline\alpha}$ is used for the product of scattering
amplitudes along the path
\[
A_{\alpha\mathbf{p}
\overline\alpha}=S_{\alpha,\alpha_n}S_{\alpha_{n},\alpha_{n-1}}
\cdots S_{\alpha_1,\overline\alpha}\,.
\]
One can now make a sequence of reductions \cite{Wilson,BHW}:
The unipotency property, $\sigma^2=1$, of the scattering matrix
(see Theorem \ref{T:ssquared}) leads to
\[
\sum_\alpha
S_{\alpha,\alpha_n}S_{\alpha_{n},\alpha_{n-1}}
\cdots S_{\alpha_1,\overline\alpha}
=\delta_{\alpha_1,\overline\alpha_n}S_{\alpha_{n},\alpha_{n-1}}
\cdots S_{\alpha_2,\alpha_1}\,.
\]
This  allows a massive inductive reduction of the sum
(\ref{E:sum}), in the course of which one must also combine
 the $\alpha$-dependent factors
$\int^{l_\mathbf{p}+2L_\alpha}_{l_\mathbf{p}}K_0(t,x)\,dx$.
 One eventually arrives at
 the following representation of (\ref{E:sum}):
\begin{equation}
\frac12 \sum_{\alpha}S_{\alpha\overline\alpha} \int_0^\infty
K_0(t,x)\,\rmd x. \label{bounces}
\end{equation}
From \eref{bounces} and \eref{E:heat} one obtains \eref{strace}
as the total contribution of all these ``bounce'' paths. This
finishes the proof of  Theorem  \ref{Ssum}.
\end{proof}

{\it Example.}
 It is well known (e.g., \cite{KSm}) that
 at a Kirchhoff vertex of degree $d_v$ the scattering matrix is
\begin{equation}
\sigma^{(v)}_{ef}= \frac2{d_v} - \delta_{ef}\,.
 \label{krefl}\end{equation}
For an entirely Kirchhoff  graph, therefore, we have
\begin{equation}
\frac14\sum_\alpha S_{\alpha\overline\alpha} =
 \frac14 \sum_{v=1}^V \left(\frac2{d_v}-1\right)d_v = \frac12\,V -
\frac12\,E, \label{kscatind}\end{equation} since every edge is
incident on two vertices. This reproduces Roth's formula
\eref{Ktrace}. By a similar calculation, or by appealing to
Proposition \ref{L:odd}, one sees that in the case of a graph all of
whose vertices are of the anti-Kirchhoff type, the resulting term is
the negative of~\eref{kscatind}.

  \begin{proposition}\label{constterm}
 In the context of Proposition \ref{H_ind},
  $\Tr K - \Tr K^\prime$ is equal to twice the constant term in the
  small-$t$ expansion of $\Tr K$.
\end{proposition}

 \begin{proof}
 Proposition \ref{L:odd} and Theorem \ref{Ssum} imply
 that the constant term in $\Tr K$ is the negative of that
in $\Tr K^\prime$.
 The proof of  Proposition \ref{H_ind} shows that only these
constant terms  survive when the traces are subtracted.
\end{proof}

\begin{corollary}\label{C:index}
Under the conditions of Theorem \ref{Ssum} and Proposition \ref{H_ind},
  the following
alternative representations hold for the index of~$A$:
\begin{equation}
\ind A = \frac12 \sum_{\alpha}S_{\alpha\overline{\alpha}}=
 \Tr K - \Tr K'
 =E-\sum_v\dim P_v=E-p,
 \label{E:altern}
\end{equation}
where $H=A^*A$ is the factorization of the Hamiltonian in accordance
with (vii) in Theorem \ref{T:ssquared}, $K$ and $K'$ are the heat
kernels of $H$ and $A^*A$, $E$ is the number of undirected edges in
$\G$, $P_v$ is the projector onto the Dirichlet part of the vertex
conditions, and $p$ is the total number of vertex conditions not
involving derivatives. Furthermore, $\Tr K - \Tr K'$ can be
read off from the asymptotics of a single heat kernel by virtue
of  Proposition~\ref{constterm}.
\end{corollary}
\begin{proof}
The last two equalities are quoted from \eref{projindthm} for
completeness.
 The rest of the corollary summarizes the results of this
subsection.
 \end{proof}

\subsection{The Euler characteristic} \label{ssec:Euler}

We  now look at the special situation of Kirchhoff conditions to
see the implications of Corollaries \ref{C:index}
 and \ref{C:index1} there.

It is convenient to consider first-order operators
on quantum graphs as defined in
terms of differential forms, rather than functions,
 introducing thus an analogue of the de Rham complex.
Therefore, we henceforth identify $A$ with
$\,\rmd\colon\Lambda^0(\Gamma) \to \Lambda^1(\Gamma)$ and $A^*$ with
$\delta\colon\Lambda^1(\Gamma) \to \Lambda^0(\Gamma)$ (on
appropriately restricted domains).

\begin{theorem} \label{Gencase}
Let $H_K$ and $H_A$ be the Kirchhoff and anti-Kirchhoff Laplacians
on a compact quantum graph $\Gamma$ (acting in $\Lambda^0(\Gamma)$
and $\Lambda^1(\Gamma)$ respectively), $K_K$ and $K_A$ be the
corresponding heat kernels, and $\,\rmd$ and $\delta= \rmd^*$ be the
external derivative operators (with the domains defined by the
continuity conditions for $\rmd$ and the sum of values equal to zero
at each vertex for $\delta$), so that
\begin{equation}
H_K = \delta \rmd \quad \mbox{and} \quad H_A= \rmd\delta.
\label{dfactor}
\end{equation}
Then
\begin{equation} \label{indextheorem}
\ind \rmd =\Tr K_K - \Tr K_A
= V-E,
\end{equation}
the Euler characteristic of $\Gamma$.
\end{theorem}

 \begin{proof}
Either the first formula
 (with \eref{kscatind}) or the last formula
  in \eref{E:altern} can be used to
calculate the index.
 \end{proof}

\begin{corollary} \label{dimkeranti}
Let $C$ be the number of
connected components of $\Gamma$.
Then
\begin{equation}
\dim\ker\,\rmd = \dim\ker H_K = C
\label{dimkerK}\end{equation}
and
\begin{equation}
\dim\ker\delta = \dim\ker H_A = E-V+C.
\label{dimkerA}\end{equation}
In particular, in the connected case
\begin{equation}
\dim\ker\delta = E-V+1=r,
\label{fundgrp}
\end{equation}
the rank of the fundamental group of $\Gamma$.
\end{corollary}

\begin{proof} Equality \eref{dimkerK} is immediate, since
the zero modes of $H_K$ are constant on each component. Now the
index theorem \eref{indextheorem} yields \eref{dimkerA}.
\end{proof}

\begin{remark} \label{deltamodes} \rm
\Eref{dimkerA} counts the locally constant differential $1$-forms
satisfying \eref{antidir}. It is $0$ for a tree graph (where
\eref{antidir} must be violated at the leaves) and $1$ for a
cycle. In general, it counts the independent cycles in the graph.
\end{remark}

\section{Relation to the secular determinant} \label{sec:secular}

Kottos and Smilansky \cite{KSm} derived their trace formula for
the density of states of a (compact) quantum graph from a certain
secular equation,
\begin{equation}
 f(k)\equiv \det [U(k)-1]=0.
\label{seceq}\end{equation}
In terms of frequency the  density of states is
\begin{equation}
\rho(k) = \sum_{n=0}^\infty \delta(k-|k_n|)
\equiv N_0 \delta(k) + \rho_1(k),
\label{denstat}\end{equation}
where $N_0$ is the  true \emph{spectral multiplicity} of $0$ as an
eigenvalue, and
$\rho_1(k)$ is the contribution of the positive eigenvalues,
$\lambda_n=k_n{}\!^2$,  appearing with multiplicity.
(We assume there are no negative eigenvalues.
According to \cite{KSm,KN,GS,Kurasov}, for strictly positive
eigenvalues the spectral multiplicity is equal to the
multiplicity as a root of \eref{seceq}, but the situation for
$k=0$ is quite different.)
The secular function $f$ defines a distribution $\tilde\rho(k)$ on the
entire real line by
\begin{eqnarray}
 \tilde\rho[\phi]
&\equiv& \lim_{\epsilon\downarrow0} \frac1{4\pi i} \int_{-\infty}^\infty
\left( {f'(k-i\epsilon)\over f(k-i\epsilon)} - {f'(k+i\epsilon)\over
f(k+i\epsilon)} \right)
\phi(k)\, dk \\
&=& \frac12 \sum_{n=-\infty}^\infty \phi(k_n),
\label{cauchy}\end{eqnarray} the sum being over all zeros of $f$
(with multiplicity), including the possible one at $k=0$, whose
\emph{algebraic multiplicity} as a root of \eref{seceq} we shall
denote~$\tilde N$. If $k$ is a nonzero root of~$f$, then so is~$-k$.
The positive eigenvalues  are the squares of the  nonzero
roots of $f$; thus every  positive eigenvalue appears in
the sum
\eref{cauchy} twice (times its multiplicity), and hence
\eref{denstat} is correctly
reproduced on the positive axis by~$\tilde\rho$. The spectral
multiplicity $N_0\,$, however, is generally  equal neither to
$\tilde N$ nor to $\frac12 \tilde N$. Furthermore, $\tilde N$ itself
has proved difficult to calculate reliably \cite{KSm,KN,GS,Kurasov}.

From the foregoing definitions and discussion it follows that
\begin{eqnarray}
 \rho(k) &=& \cases{
\bigl( N_0-{\textstyle\frac12} \tilde N \bigr) \delta(k)
+\tilde\rho(k) & for $k\ge
0$, \cr
0 & for $k<0$ \cr}
\nonumber \\
&\equiv& h(k)\bigl[\bigl( N_0-{\textstyle\frac12} \tilde N \bigr)
\delta(k)
+\tilde\rho(k)\bigl],
\label{rhocases}\end{eqnarray}
where $h$ is the unit step function and
(because the product $h\delta$ would otherwise be ambiguous) we
stipulate that
$h\delta[\phi] \equiv \phi(0)$ (not $\frac12\phi(0)$).
In terms of the density of states,
\eref{roth} can be rewritten as
\[\int_{-\infty}^\infty e^{-k^2t} \rho(k)\,dk
= \Tr K = K_1 + K_2+K_3\,. \]
The derivation of the trace formula \cite{KSm,KN,GS,Kurasov}
makes clear that $\tilde\rho$
yields precisely the leading Weyl term and the oscillatory
(periodic-orbit) terms in  the density of states.
It follows (cf.\ discussion following \eref{roth}) that the
contribution of $\tilde\rho$ to the heat kernel is
 the  terms $K_1$ and $K_2\,$.
 The remaining term
in \eref{rhocases}, proportional to $\delta(k)$,
must therefore be responsible for precisely
the constant term $K_3$ in the heat kernel.
We reemphasize
 that the coefficient of those terms is not $N_0\,$,
the coefficient of $\delta(k)$ in the spectral density~$\rho$.
That is not a paradox:
The distributionally convergent periodic-orbit sum
$\tilde\rho$
contains another
$\delta(k)$ contribution to $\rho$ that restores consistency
with~\eref{denstat}.

On the other hand, we now know from previous sections that
$K_3\,$,
for a scale-invariant Laplacian,  is equal to half the
index of the associated first-order operator,~$A$.
The index, in turn, equals $E-p$, where $E$ is the number of (undirected)
edges and $p$ is the number of Dirichlet conditions.
Furthermore, we have $\ind A = N_0-N_0^*$,
where $N_0^*$ is the nullity of $A^*$ and hence of the operator
dual
to~$H$.
Therefore, we immediately get two  interesting identities:

\begin{corollary}\label{algmult}
Let $N_0$ and $\tilde N$ be the spectral and algebraic multiplicities of
$k=0$ for  a scale-invariant graph Laplacian, $A^*A$, and let
$N_0^*$ be the
spectral multiplicity for the dual Laplacian, $AA^*$.  Then
\begin{equation}
 \tilde  N = 2N_0 - \ind A
= 2N_0 - E +p
\label{algmult1}\end{equation}
and
\begin{equation}
\tilde N= N_0+N_0^*\,.
\label{algmult2}\end{equation}
\end{corollary}

\emph{Example 1:}  For $H_K\,$, the Laplacian of a connected Kirchhoff
graph, one has $N_0 =1$ and $\,\ind A=\ind\rmd = V-E$.
Therefore,
$\tilde N = 2-V+E$,
in agreement with \cite[corrigendum]{KN} and \cite{Kurasov}.

\emph{Example 2:}  For the pure Neumann Laplacian $H_N$ of
\sref{sec:count}, one has $N_0=E$ and $p=0$, so $\tilde N=E$.
This is correct, because $0$ appears as a root of $f$ once for each
disconnected Neumann edge.

Kurasov \cite{Kurasov} gives a convincing direct calculation of $\tilde N$
for the Kirchhoff case.  On that basis he deduces that the Euler
characteristic \eref{Ktrace} is determined by the spectrum of~$H_K\,$.
Thus the direction of the logic in \cite{Kurasov} is roughly the reverse
of that in the present paper.
Our derivation of Corollary \ref{algmult} is simpler (as well as more
general).

\section{Conclusions and additional remarks} \label{sec:concl}

We have demonstrated that the ``topological'' term in the
heat-kernel expansion for a Laplacian on a quantum graph does
indeed have an index interpretation, if the Laplacian is of the
scale-invariant class (i.e., the boundary conditions do not mix
function values and derivatives, and hence the scattering matrix
is independent of~$k$). Such a Laplacian factors into two
first-order operators, $A$ and $A^*$, defined
on domains determined by those boundary conditions.

One can calculate the index in three ways:
(1)~in the usual way, by subtracting the heat kernel of
$AA^*$ from that of $A^*A$
(Proposition \ref{H_ind} and Theorem \ref{Gencase});
(2)~by inspection of just the heat kernel of $A^*A$
(Proposition \ref{constterm} and Corollary \ref{C:index});
(3)~just by counting the number of Dirichlet-type conditions
(Theorem \ref{genorder} and Corollary \ref{C:index1}).

A general index formula has been  derived for an arbitrary
(elliptic)
differential operator on a quantum graph (Theorem \ref{genorder}).

Along the way, we have provided some properties of the
Hamiltonian
that are equivalent to the scale invariance of the vertex
conditions (Theorem \ref{T:ssquared}).

Finally, we have determined the algebraic multiplicity of $0$ as a
root of the secular equation of a generic scale-invariant  graph
Laplacian in a novel way (Corollary \ref{algmult}).

We now add a few final remarks concerning the results of the
paper:

\begin{itemize}
\item There is an elementary sense in which  the integer $V-E$
encountered in \eref{Ktrace} is associated
with an operator index.
A graph as
a purely combinatorial object (with no lengths assigned to the
edges) is described in graph theory by the
\emph{incidence matrix}, whose rows are indexed by the vertices
and its columns by the edges. Each matrix entry is  equal
to either $0$, $1$, or~$2$, depending on whether that edge is not or
is  attached to that vertex or forms a loop there.
 Then it is easy to see that (just because of the matrix's
dimensions) the index of the incidence matrix is equal to $E-V$.

\item
The restriction to ``local'' vertex conditions in Theorem
\ref{T:Sa}
 and elsewhere has very little content.
Indeed, the structure of the graph enters the problem only
 through the vertex conditions, and one could \emph{define} a
vertex as a subset of edges that are related by such conditions.
Alternatively, one could think that all
vertices of a quantum graph have collapsed into a single one
(creating a ``rosette'' consisting of  one vertex and $E$ cycles
attached).
Since all (whether previously local or nonlocal) vertex
conditions refer to  this single vertex,  all conditions  have
become local.
This is impossible only if we need to enforce some specific type
of vertex conditions, e.g., the Kirchhoff ones;
after the graph collapses to a rosette, the vertex conditions
will generally no longer be of that type.

\item  Kirchhoff conditions on $\Lambda^0(\Gamma)$ and
anti-Kirchhoff conditions on $\Lambda^1(\Gamma)$ seem to be quite
natural, whereas the interchanged conditions look rather
unnatural. The situation for manifolds is different: there one can
either impose Neumann conditions on forms of even degree and
Dirichlet conditions on forms of odd degree, or do the reverse.
The two choices correspond to two different cohomology theories
for the manifold, ``absolute'' and ``relative''
\cite{Gilkeybook,Gilkey}.

\item Carlson \cite{Carlson} also constructs second-order
self-adjoint operators on quantum graphs in terms of first-order
operators, but his construction is rather different from ours. More
pertinent are the remarks of Friedman and Tillich \cite{FT} and
Exner and Post \cite{EP} that derivatives on quantum graphs should
be treated as 1-forms (or vector fields).

\item
An extension of the first-order formulation of the index theorem
 to operators $H$ with a
nontrivial Robin part, $C_v\ne 0$, is not to be expected.
 As  pointed out by V.~Kostrykin, in general either
$H$ or its dual will have negative spectrum, whereas operators of
form $A^*A$ and $AA^*$ must both be positive.
 This impossibility of factorization is also contained in the
statement (vii) of Theorem~\ref{T:ssquared}. On the other hand,
inserting ``Robin'' operators $\Lambda_v$ into the definition of
$H$ does not change the $t$-independent terms of $\Tr K_H$ and
its dual, which are the ingredients of the formula for the index
of the second-order operator.

\end{itemize}

\ack We thank P. Gilkey, J. Harrison, J. Keating, K.~Kirsten,
P.~Kurasov, M.~Pivarski, R.~Schrader, A.~Terras, B.~Winn, and
especially G.~Berkolaiko, J.~Bolte, and V.~Kostrykin for various
forms of advice or assistance. This research was supported in part
by NSF grants PHY-0554849, DMS-0406022, DMS-0648786 and the
program Analysis on Graphs and Their Applications at the Isaac
Newton Institute for Mathematical Sciences at Cambridge
University.

\end{document}